\pdfoutput=1
\def\arxivversion{}
\documentclass[12pt]{article}

\usepackage[utf8]{inputenc}
\usepackage[T1]{fontenc}
\usepackage{lmodern}
\usepackage[margin=1in]{geometry}
\usepackage{amsmath,amssymb}
\usepackage{graphicx}
\usepackage{csquotes}
\usepackage[backend=biber,style=apa,natbib=true]{biblatex}
\usepackage{hyperref}
\usepackage{orcidlink}

\title{An Ontology-Based Approach to Optimizing Geometry Problem Sets for Skill
Development}

\ifdefined\anonymous
\author{}
\else
\author{%
  Michael Bouzinier\,\orcidlink{0000-0002-3161-5601}\thanks{Harvard University, MA}%
  \thanks{Forome Association, MA}%
  \thanks{IDEXX Laboratories, ME}%
  \thanks{Russian School of Mathematics, MA} \and
  Sergey Trifonov\footnotemark[2] \and
  Matthew Chen\,\orcidlink{0009-0007-2272-6215}\thanks{Lexington High School, MA} \and
  Tarun Venkatesh\,\orcidlink{0009-0002-6148-2982}\thanks{Park Hill South High School, MO} \and
  Lielle Rifkin\,\orcidlink{0009-0006-0169-2944}\footnotemark[4]\thanks{Brookline High School, MA}%
}
\fi
\date{}

\begin{document}

\maketitle
\thispagestyle{empty}

\begin{abstract}
Euclidean geometry has historically played a central role in cultivating
logical reasoning and abstract thinking within mathematics education, but has
experienced waning emphasis in recent curricula. The resurgence of interest,
driven by advances in artificial intelligence and educational technology, has
highlighted geometry's potential to develop essential cognitive skills and
inspired new approaches to automated problem solving and proof verification.
This article presents an ontology-based framework for annotating and optimizing
geometry problem sets, originally developed in the 1990s. The ontology
systematically classifies geometric problems, solutions, and associated skills
into interlinked facts, objects, and methods, supporting granular tracking of
student abilities and facilitating curriculum design. The core concept of
`solution graphs'---directed acyclic graphs encoding multiple solution pathways
and skill dependencies---enables alignment of problem selection with
instructional objectives.
The framework has been tested in practice through the annotation of thousands
of problems over three decades. We contend that our approach addresses
longstanding challenges in representing dynamic, procedurally complex
mathematical knowledge. We conclude by articulating a research agenda: the open
problems of automated problem annotation and solution validation, whose
resolution would reduce the time teachers spend validating student work and
enable interactive feedback for self-learners.
\end{abstract}


\ifdefined\arxivversion
\begin{center}
  \includegraphics[width=0.28\linewidth]{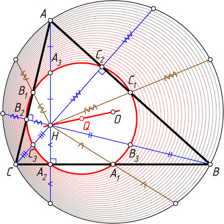}
\end{center}
\fi

\section{Introduction}

\subsection{Background and Motivation}

Euclidean geometry has long been important to mathematical education, developing
students' logical reasoning abilities and familiarity with abstract concepts
through intuitive, visual problem-solving contexts. Geometry's foundational role
in mathematics was cemented by Hilbert \parencite{hilbert1971foundations}, and
its significance for mathematical thinking has been emphasized by Freudenthal
\parencite{freudenthal1973mathematics,freudenthal1983didactical}. A
comprehensive review of recent approaches to geometry education is provided by
Sinclair et al.\ \parencite{sinclair2017recent}.

However, throughout the latter half of the twentieth century and into the early
twenty-first century, geometry experienced a decline in prominence within many
Western mathematics curricula. This trend may be attributed, in part, to the
perception that geometry had limited practical application in engineering and
technology, especially as computing proliferated and reshaped technical tasks.
Another contributing factor may have been the relative difficulty of effectively
integrating geometry into computer-aided educational techniques during the
initial rise of such technologies, at least until the 2020s.

Today, both rationales appear increasingly outdated. Advancements in AI-based
software have begun to automate routine and even intuitive tasks, elevating the
value of human skills such as abstract reasoning and logical
problem-solving---skills perfectly developed through geometric education. At the
same time, given persistent student struggles with proof construction
\parencite{stylianides2024proof}, the widespread omission of reasoning-and-proving
opportunities from textbooks \parencite{thompson2012opportunities}, and the fact
that many students reach proof-oriented geometry without having attained even the
informal-deduction level of geometric thinking \parencite{vojkuvkova2012role},
effective proof education is more vital than ever. Meanwhile, significant progress in educational
technology now makes dynamic visual tools, automated proof verification, and
interactive geometrical construction widely available, thereby enhancing both
accessibility and effectiveness of geometry education.

Despite broad advances in AI-driven mathematical reasoning
\parencite{li2024survey,murphy2024autoformalizing,polu2020generative}, geometry
continues to pose unique challenges for automation, motivating the development of
robust, machine-readable ontologies to support both learners and AI systems.

In light of this renewed importance, we present a comprehensive ontology and
methodology for geometry problem annotation, originally developed in the early
1990s for educational software \parencite{donskaya2018reshat}. While the core
framework was established over three decades ago, we argue that it has renewed
relevance in the context of modern educational technology and artificial
intelligence applications. Recent developments in automated theorem proving and
AI-assisted mathematics education, exemplified by systems like AlphaGeometry
\parencite{trinh2024solving}, suggest that structured approaches to problem
annotation and solution validation are becoming increasingly important.

Our contribution is threefold: (1) we document a mature ontological framework for
geometry problem classification that has been tested with thousands of problems
over multiple decades; (2) we demonstrate how solution graphs can systematically
support curriculum design and skill development tracking; and (3) we articulate
the open problems for automated solution validation and interactive learning
support.

\subsection{Historical Context}

Before presenting our contribution, it is important to review key historical
developments in both geometry education and the teaching of logical reasoning.
Euclid is often credited with introducing the axiomatic method that underpins
much of modern mathematical proof. In the \emph{Elements}
\parencite{euclid1956elements}, he set out definitions of basic terms such as
points, lines, and planes, together with postulates and common notions, and used
deductive reasoning to develop geometric statements systematically. Building on this foundation, Hilbert later refined the
axiomatic approach, organizing the axioms of geometry into five
groups---connection, order, parallels, congruence, and continuity---thereby
establishing a more rigorous, formal structure for geometric proof
\parencite{hilbert1971foundations}.

Throughout the twentieth century, mathematics education underwent significant
changes. Most notably, the ``New Math'' movement of the 1960s, influenced by the
formalist and set-theoretic perspectives of the Bourbaki group, sought to
restructure curricula to emphasize rigor, abstraction, and theorem proving
\parencite{kilpatrick2012newmath,kline1973why}. Later educational theorists, such
as Raymond Duval, advocated for a semiotic (sign-based) approach to mathematics
instruction that takes students' cognitive processes into account
\parencite{duval2006cognitive}. Freudenthal, meanwhile, saw proof as part of
the broader activity of mathematizing---organizing reality ``with mathematical
means'' \parencite{freudenthal1973mathematics,freudenthal1983didactical}. In contemporary pedagogy, the
trend increasingly favors inquiry- and proof-based strategies, where students are
encouraged to construct, critique, and explain arguments, rather than merely
reproduce formal proofs
\parencite{stylianidesaj2007proof}.

A particularly influential development took place in the Soviet Union, where N.N.
Konstantinov and colleagues
\parencite{gerver1965zadachi,konstantinov1971preprint} pioneered an approach in
which every mathematical fact not designated as an axiom---and therefore subject
to proof---should be presented to students as a problem to be solved.
Accordingly, students were tasked with proving all theorems in the curriculum
themselves, rather than simply reading and reciting pre-written proofs. The
curriculum thus became a sequenced collection of problems, with students guided
by teachers, junior teachers, and teaching assistants as they constructed their
own proofs. Conventional, non-theorem problems were also included in the
curriculum after students had established the foundational results through proof.

Meanwhile, theoretical advances in Western Europe were contributing additional
perspectives on how students develop geometric understanding. The van Hiele model
\parencite{vanhiele1959problem,vanhiele1986structure,gutierrez1991alternative},
in particular, provided a
framework for categorizing levels of geometric reasoning and emphasized the
importance of both developmental stages and inquiry-based pedagogy. Although the
van Hiele model is now recognized internationally, it was in many respects
adopted more thoroughly in Soviet mathematics curricula than in those of many
Western countries
\parencite{kolmogorov1977geometry6,kolmogorov1977geometry7,kolmogorov1977geometry8,usiskin1982van}.

The unique ecosystem of mathematics education that emerged in the Soviet Union
was shaped by the combined influences of Hilbert's formalism, the Bourbaki
group's structural approach, and the van Hiele model's focus on cognitive
development. Within this environment, our team developed a software tool to
support high school teachers in teaching Euclidean Geometry
\parencite{donskaya2018reshat}. The tool enabled educators to create problem sets
that demonstrated the practical application of specific problem-solving
techniques and provided students with supplementary assignments to reinforce
these skills. Initially released as a standalone DOS application, the software
included approximately 5,000 geometry problems in Russian, many of which were
also included in the book by Sharygin and Gordin
\parencite{sharygin2001sbornik}. In subsequent years, a new development team
transitioned the tool to a web-based platform, expanding the problem database to
about 16,000 problems \parencite{mccme2026zadachi}, with many published in the
book by Gordin \parencite{gordin2006geometriya}.

To optimize the selection of problems for constructing these sets, we designed a
comprehensive ontology for problem annotation. This ontology includes roughly two
dozen annotation classes and subclasses, covering attributes such as difficulty,
problem type, topic, purpose, and more. Further details about the ontology are
provided in the Supplemental Note 1.

\subsection{Ontologies in Euclidean Geometry}

Current ontologies in Euclidean geometry represent elements such as points,
lines, angles, and figures, along with their properties and logical relationships
between them in hierarchical structures in a formal and machine-readable way
\parencite{ontomath2014}. They are usually expressed in formal
knowledge representation languages such as OWL or RDF, which enable the formal
description of categories, subcategories, and relationships between objects using
logical rules, allowing computers to interpret and reason about these
connections. Many modern ontologies also incorporate links to
knowledge graphs, which enable relationships to be represented as interconnected
nodes and edges that include both definitions and dependencies among concepts
\parencite{ontomath2014}.

Despite advances, several issues persist in current ontological approaches to
Euclidean geometry. One major challenge is the limited representation of
procedural and inferential knowledge. While many ontologies effectively capture
the fixed structure of geometric entities and relationships, they often fail to
model the dynamic reasoning involved in solving or proving geometric problems
\parencite{lange2013formalization}. As a result, these systems can represent what
geometric objects are, but not how they interact within logical deduction or
construction. Another issue is the trade-off between expressivity and
computational efficiency. Ontologies that are expressive enough to encode complex
geometric relations, such as congruence or parallelism, typically become
computationally demanding or incompatible with automated reasoning systems
\parencite{elizarov2022ontomath}.
Moreover, many ontological frameworks for
geometry experience rigidity and scalability issues. Adding new configurations or
extending the ontology to related domains often requires extensive
restructuring, limiting reuse and adaptability.
The integration of different types of representations, such as written problem
statements, mathematical symbols, and geometric diagrams, remains a major
challenge because bringing these varied forms of data together in one ontological
model requires keeping their meaning consistent across all formats.
Finally, educational applications of mathematical
ontologies---geometry ontologies in particular---remain underdeveloped, with few
models capable of supporting feedback loops or automated reasoning about
students' intermediate steps in problem solving
\parencite{lalingkar2014ontology,stancin2020ontologies}.

\section{Methodology}

We now introduce our ontology-based framework, building on the historical and
technological foundations discussed above.

\subsection{Ontology Designed to Annotate Problem Solutions}

The annotation of geometry problems and their solutions is most effective when
supported by a rigorous, well-structured framework. To ensure that educational
content can be consistently categorized, efficiently searched, and readily
adapted to evolving curricula, it is essential to establish a formal
ontology---a systematic classification of all relevant elements and their
relationships within the domain. We have developed such an ontology specifically
for the context of student problem solving, identifying and organizing the key
components involved in the process. These elements are grouped into three main
classes:

\begin{itemize}
  \item \textbf{Facts}: Euclidean axioms, derived theorems, lemmas, and other
    provable statements within the Euclidean geometry framework.
  \item \textbf{Geometric Objects}: Figures and concepts that are either given in
    problem statements or constructed during the problem-solving process.
  \item \textbf{Methods}: Specific techniques or strategies that can be applied
    to arrive at a solution.
\end{itemize}

Within the triadic structure of our ontology---facts, objects, and
methods---the designation of \textbf{objects} as a distinct category plays a
fundamental role. The composition of the object category is typically the most
stable aspect of the ontology. While the lists of facts (such as theorems or
lemmas) and methods (problem-solving strategies) gradually stabilize with ongoing
development and extensive annotation of problem collections, they may remain
somewhat flexible and subject to extension and revision as the ontology evolves.
In contrast, changes to the set of objects are exceedingly rare, yet when they do
occur, they can significantly alter the overall structure of the subject domain.
Objects serve as the conceptual building blocks; modifying this set has deep
ramifications for the organization and interpretation of both facts and methods
across the ontology. For instance, introducing a new geometric object---such as a
``nine-point circle''---requires re-examining related facts and may enable new
methods, illustrating how fundamental changes to the object category cascade
through the ontology.

Notably, while ``Methods'' served historically as the core building blocks of the
ontology, there is a nuanced difference between a \textbf{method} as an element
defined in the ontology and a \textbf{skill} as used in education. Methods refer
to well-defined, specific techniques---such as constructing an altitude, using
the area method, or applying a sequence of steps to prove congruence. Skills, in
the pedagogical sense, encompass more than correct execution of methods; they
include the ability to recognize when and how to apply these methods in varied
contexts, to recall and deploy relevant facts, to identify pertinent objects, and
to creatively combine multiple elements in problem-solving. In short, a student's
skill is demonstrated in their mastery of selecting, integrating, and applying
methods, facts, and objects to solve a wide array of problems.

This ontology also served as a foundation for the book \emph{Every Student Must
Know This} \parencite{gordin2003eto}.

\subsection{Problem Solutions Graph Representation}

It is important to acknowledge that many problems can be approached and solved
using various methods. Good educators should encourage diverse approaches to
problem solving, rather than insist on solutions that mirror those demonstrated
in class. Nevertheless, curricula seek to develop specific skills, so problem
selection should align with these learning objectives.

Throughout the curriculum, students continually reinforce skills they already
possess while being introduced to new ones. For effective skill illustration or
practice, problems should meet these criteria:

\begin{enumerate}
  \item \textbf{Existing Skills}: The problem should require only those skills
    students have already acquired.
  \item \textbf{Target Skill Integration}: The skill being taught should play an
    essential role in the solution.
  \item \textbf{Skill Necessity or Efficiency}: There should be no reasonable
    alternative solution that excludes the target skill, or else the skill should
    significantly simplify the solution.
\end{enumerate}

\begin{figure}[htbp]
  \centering
  \includegraphics[height=0.82\textheight]{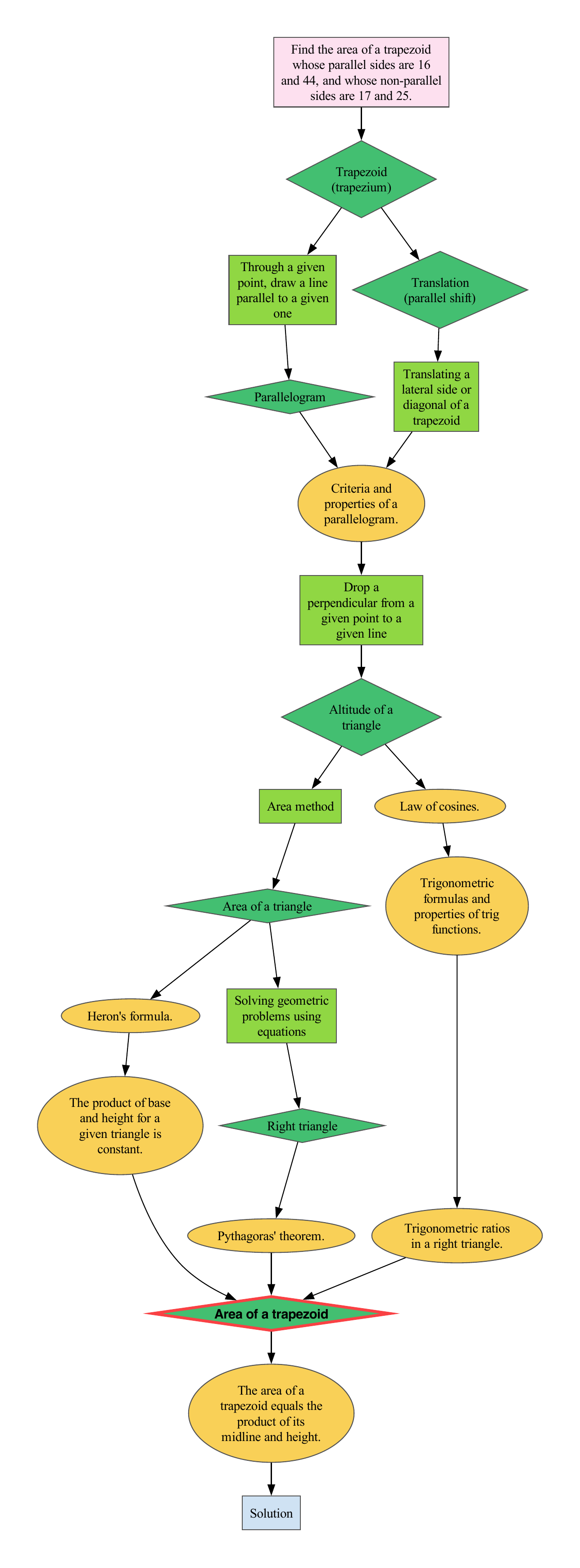}
  \caption{Illustration of a solution graph: a directed acyclic graph of the
  facts, objects, and methods that connect a problem statement to its solution.}
  \label{fig:solution-graph}
\end{figure}

To support efficient problem selection and curriculum alignment, we represent
solutions as \textbf{Solution Graphs}---directed acyclic graphs (DAGs) whose
nodes correspond to skills (facts, objects, or methods). Effective problem design
requires a path through the graph that passes through required skills and does
not allow shortcuts that bypass newly introduced skills (see Figure~\ref{fig:solution-graph}
for illustration). This approach enables teachers to track the learning
trajectory and ensure curricular coherence at every stage.

Recent AI systems for geometry, such as AlphaGeometry
\parencite{trinh2024solving}, similarly represent solutions as proof graphs for
automated search and validation.

Our experience over the past thirty years has shown that, in some cases, this
methodology empowered self-motivated students to develop significant
problem-solving skills independently, with only occasional teacher consultation.
The granularity of method-based tracking and automated feedback made
self-directed skill development and mastery possible.

\section{Discussion}

\subsection{Today's Relevance: from building sets to validating proofs}

The original tool, though innovative, remains a niche product. This is primarily
due to two factors. First, constructing problem sets is typically a task for
curriculum developers, not an everyday responsibility for teachers. Teachers
might be creating personalized training assignments, but this happens
infrequently in practice. Secondly, Euclidean Geometry has become less prominent
in modern Western education. Ironically, the latter factor also underscores why
the tool's foundation might hold greater relevance today than it did 30 years ago.

Euclidean Geometry teaches students both logical reasoning and manipulation with
abstract concepts. At the same time, two-dimensional geometry is intuitive and
allows for easy illustration of both logical frameworks and abstractions through
text and geometric drawings. This makes it an ideal medium for teaching these
skills to young students when their brains are most receptive. In an era where AI
is expected to offload many routine tasks, proficiency in logical reasoning and
high-level abstraction is becoming increasingly crucial. Thus, Euclidean Geometry
presents itself as a unique and valuable subject in contemporary school
curricula. Recent developments in AI, such as AlphaGeometry
\parencite{trinh2024solving}, underscore both the feasibility and the pedagogical
potential of automated geometry problem solving and proof generation.

One barrier to its popularity might be the challenge of verifying solutions to
geometry problems, a task requiring significant expertise and tedious work. These
problems often require proofs that can be composed in multiple ways. Demanding
that young students produce strictly formalized proofs for automated validation
is impractical. Moreover, geometric construction problems, usually described less
formally than proofs, further complicate automatic validation.

We might speculate that a software tool to assist teachers in validating student
solutions would make popularizing the teaching of Euclidean Geometry more
feasible, leading to including the subject in the majority of school curricula.
We believe the methodology for problem solution annotation that we developed
three decades ago could significantly aid in automating solution verification.

With the explosion in popularity of AI in general, and large language models
(LLMs) in particular, many recent works have focused on generating and validating
formal proofs for mathematical problems, including Euclidean geometry. For
instance, AlphaGeometry \parencite{trinh2024solving} underscores the critical role
of structured geometric knowledge in enabling automated proof search and
validation. Their system relies on a comprehensive database of geometric
theorems, facts, and construction techniques, and effectively organizes proof
steps as a DAG. However, the authors stop short of formalizing their technique as
an ontology-based approach; that is, the explicit classification and
systematization of geometric knowledge into foundational objects, facts, and
methods remains implicit rather than formally defined. Furthermore, in
AlphaGeometry, the primary purpose of constructing a DAG is to search for a single
valid proof, rather than to classify and organize all possible approaches to
solving a problem, as emphasized in our ontology-based methodology. As noted
by Sinclair et al.\ \parencite{sinclair2017recent}, digital technologies have
become mainstream in geometry education, with new approaches emerging in
assessment and feedback.

\subsection{Future Directions}

We foresee two milestones on the path toward automating solution validation.
The first is computer-aided validation that checks routine student proofs and
constructions, reducing the time teachers spend on grading while leaving
unusual or creative solutions to human review. The second is interactive
feedback for self-learners that flags incorrect or irrelevant statements as
they write and offers targeted hints. Existing interactive theorem provers
such as Coq and Lean can verify fully formal proofs, but current AI systems
still struggle to formalize typical classroom geometry solutions from
informal text and diagrams, with only modest success rates
\parencite{li2024survey}. Our ontology-based representation may narrow this
gap by treating formalization as semantic parsing, a task where large
language models already perform well.

Because our approach is based on representing each problem's solution as a
graph of facts, objects, and methods, further progress toward these
milestones will depend on advances in three open problems.
First,
solution graphs need to be linked to existing formal proof tools so that
they can check individual steps or whole proofs. Second, problem annotation
needs to be automated by constructing solution graphs from problem
statements and diagrams. Third, a student's solution must be mapped to a
known solution graph; like annotation, this is a semantic parsing task for
which current LLMs are well suited. In U.S. curricula, proofs are often
written in a statement--reason form rather than free text, which simplifies
parsing, whereas construction problems lack a standard format and are
harder to process automatically. A statement--reason proof, viewed as a
sequence of (statement, reason) pairs, can be translated into a formal
language suitable for proof checking; theorem-proving environments such as
Lean
\parencite{demoura2021lean,murphy2024autoformalizing,qian2025lean,song2025leangeo}
and Coq
\parencite{narboux2018synthetic,coq2023proof,coq2023refman}
are natural targets. If validation operates on selected subsequences of such
pairs, and uses prior knowledge of the relevant solution graph, it can label
each statement as incorrect or unproven, correct but justified by an
invalid reason, correct but probably irrelevant, or correct and relevant,
and thus decide whether the overall solution is acceptable.

Interactive feedback, especially for self-learners, will likely need
to support several ways of entering statement--reason proofs. One
option is constructed statements chosen from a controlled text
interface; this scaffolds students who struggle with mathematical
language and simplifies automated checking. Another is to link each
statement to actions on a geometric figure in dynamic geometry
software. This is intuitive for simpler problems, but for advanced or
degenerate configurations it can break down, for example when key
objects lie outside the visible part of the diagram. A third option is
free-form, write-in text, which offers the greatest flexibility and is
essential for challenging problems; recent work in semantic parsing
suggests that such statements can sometimes be translated into formal
Lean~4 code \parencite{murphy2024autoformalizing,song2025leangeo}. The
most effective mix of these input modes, and the extent to which they
can support detailed educational feedback, remains an open question.

None of these components has yet been built or evaluated; together they
define a research agenda. Our first step is the annotation problem: an
empirical study of automated annotation of problems from our corpus is in
preparation and will be reported separately.

\ifdefined\anonymous\else
\section*{Acknowledgments}

The majority of this work was performed in 1991--1996 by R.K.\ Gordin, M.\
Bouzinier, S.\ Trifonov, and I.F.\ Sharygin. Unfortunately, I.F.\ Sharygin passed
away on the 12 of March, 2004 and R.K.\ Gordin is unable to continue working on
this project because of reasons beyond our control. We acknowledge a major
intellectual contribution by I.F.\ Sharygin and the primary role that R.K.\
Gordin has played in developing the product, however for the reasons above we
cannot list them as the authors.

The original software would not have been possible without the support of Arkady
Volozh, discussions with Grigory Kondakov, Alexey Kanel-Belov, Victor Prasolov,
and many colleagues in Math Education in 1990s Moscow. More recently, rewarding
conversations with Alex Karezin and Ilya Rifkin inspired the first draft of this
new manuscript. Over the last year, many long discussions with Michael Bukatin
have shaped its final form.

\section*{Author contributions}

M.B.\ and S.T.\ contributed to problem conceptualization and tool design.
S.T.\ developed the majority of the software code, with contributions from
M.B. The initial draft of the manuscript was written by M.B. M.C.\ and T.V.\
reviewed the literature on the use of AI in mathematics education and
automated proofs, and drafted the corresponding sections. L.R.\ reviewed the
literature on methodologies for teaching Euclidean geometry and drafted the
related sections. All authors reviewed, commented on, and edited the final
manuscript. All authors approved the final version.

\section*{Disclosure statement}

The authors report there are no competing interests to declare.

\section*{Funding}

No funding was obtained for the work reported in this article.

\section*{Declaration of generative AI use}

Generative AI tools (OpenAI GPT~4.1 and Anthropic Claude~4 Sonnet) were used
in the preparation of this manuscript to check grammar, improve style, and
check the text for consistency. AI assistance was also used, under author
supervision, to convert the manuscript to \LaTeX{} and to audit the
bibliography, in the course of which every reference was verified against its
published source and corrected where necessary; the audit record is retained
by the authors. Generative AI was not used to produce the research content of
this article. The authors reviewed and take full responsibility for the
content of this article.

\section*{Data availability statement}

No datasets were generated or analysed in the preparation of this article.
The ontology vocabularies documented in the supplement (the catalogs of
geometric facts, objects, and methods, in Russian and English) are available
from the corresponding author on reasonable request; a curated public release
is in preparation.

\fi

\printbibliography

\end{document}